\documentclass{article}
\usepackage{amsmath,amssymb,amsthm}
\usepackage{graphics}

\newtheorem{theorem}{Theorem}[section]
\newtheorem{lemma}{Lemma}[section]

\newcommand{\dist}{\mathop{\rm dist}\nolimits}

\newcommand{\cg}{\mathop{\rm CG}\nolimits}
\newcommand{\ar}{\mathop{\rm area}\nolimits}
\newcommand{\ko}{\mathop{\rm K}\nolimits}
\newcommand{\iq}{\mathop{\rm IQ}\nolimits}
\newcommand{\dt}{\mathop{\rm DT}\nolimits}
\newcommand{\vd}{\mathop{\rm VD}\nolimits}

\newtheorem{cor}{Corollary}[section]

\title {Five Essays on the  Geometry of L\'aszl\'o Fejes T\'oth}

\author {Oleg R. Musin\thanks{This research is partially supported by the NSF grant DMS-1400876 and the RFBR grant 15-01-99563.}}

\begin{document}
\date{}
\maketitle

\begin{abstract} In this paper we consider the following topics related to results of L\'aszl\'o Fejes T\'oth:
(1) The Tammes problem and Fejes T\'oth's bound on circle packings;
(2) Fejes T\'oth's problem on maximizing the minimum distance between antipodal pairs of points on the sphere;
(3) Fejes T\'oth's  problem on the maximum  kissing number of packings on the sphere;
(4) The Fejes T\'oth -- Sachs problem on the one--sided kissing numbers;
(5) Fejes T\'oth's papers on the isoperimetric problem for polyhedra.
\end {abstract}

\section{Tammes' problem and  Fejes T\'oth's bound on circle packings}

\subsection{Tammes' problem}

We start with the following classical problem: How should $N$ points be distributed on a unit sphere so that the minimum distance between two points of the set attains its maximum value $d_N$? This problem was first asked by the Dutch botanist Tammes \cite{Tam} while examining the distribution of openings on the pollen grains of different flowers. This question {is} also known as the problem of the ``inimical dictators'' \cite{Mesch}, namely {\it ``where should $N$ dictators build their palaces on a planet so as to be as far away from each other as possible?''} The problem is equivalent with the problem of densest packing of congruent circles on the sphere (se {\it{e.g.}} \cite[Section 1.6: Problem 6]{BMP}): How are $N$ congruent, non-overlapping circles distributed on a sphere when the common radius of the circles has to be as large as possible? The higher dimensional analogue of the problem has applications in information theory \cite{vdW61}. This justifies the terminology that a finite subset $X$ of ${\Bbb S}^{n}$ with
$$\psi(X):=\min\limits_{x,y\in X,x\ne y}{\dist(x,y)}$$
is called a spherical $\psi(X)$-code.

Tammes' problem is presently solved only for $N\le14$ and  $N=24$. L. Fejes T\'oth \cite{FeT0} solved the problem for $N=3,4,6,12$. Sch\"utte and van der Waerden \cite{SvdW1} settled the cases $N=5,7,8,9$. The cases $N=10$ and $11$ were solved by Danzer \cite{Dan}\footnote{Actually, Danzer's paper \cite{Dan}  is the English translation of his Habilitationsschrift ``Endliche Punktmengen auf der 2-sph\"are mit m\"oglichst gro{\ss}em Minimalabstand''. Universit\"at G\"othingen, 1963.} (see also the papers by B\"or\"oczky \cite{Bor11} for $N=11$ and H\'ars \cite{Hars} for $N=10$). Robinson \cite{Rob} solved the problem for $N=24$. In my recent papers with Tarasov \cite{MT,MT14} we gave a computer--assisted solution for $N=13$ and $N=14$.

Robinson extended Fejes T\'oth's method and and gave a bound valid for all $N$ that is sharp besides the cases $N=3,4,6,12$ also for $N=24$. The solution of all other cases is based on the investigation of the so called contact graphs associated with a finite set of points.
For a finite set $X$  in ${\mathbb S}^2$ the {\it contact graph} $\cg(X)$ is the graph with vertices in $X$ and edges $(x,y), \, x,y\in X$, such that $\dist(x,y)=\psi(X)$. The concept of contact graphs was first used by Sch\"utte and van der Waerden \cite{SvdW1}. They used the method also for the solution the thirteen spheres (Newton--Gregory) problem \cite{SvdW2}.

In Chapter VI of the book \cite{FeT} the concept of irreducible contact graphs is considered in details. 
The method of irreducible spherical contact graphs was used also \cite{HabvdW,vdW,Bor,BS,BS14} for obtaining bounds for the kissing number and Tammes problem.

\subsection{The Fejes T\'oth bound}

Now we consider a theorem on bounds of equal--circle packing and covering of a sphere proved by L\'aszl\'o Fejes T\'oth in 1943 \cite{FeT0,FeT,FeTbook2}.

\begin{theorem}[L. Fejes T\'oth \cite{FeT0}] If $\Theta$ is the density of a packing of the unit sphere ${\Bbb S}^2$ in ${\Bbb R}^3$  by $N$ congruent spherical caps then
	$$
	\Theta\le \frac{N}{4}\left(2-\csc{\omega_N}\right),  \; \mbox{ where } \; \omega_N:=\frac{N\pi}{6N-12}
	$$
If $\Omega$ is the density of a covering of \, ${\Bbb S}^2$ by $N$ congruent spherical caps then
		$$
		\Omega\ge \frac{N}{2}\left(1-\frac{1}{\sqrt{3}}\cot \omega_N\right).
		$$
\end{theorem}

Denote by $A(n,\varphi)$ the maximum cardinality of a $\varphi$--code in ${\Bbb S}^{n-1}$.  In other words, $A(n,\varphi)$ is the maximum cardinality of a packing in ${\Bbb S}^{n-1}$ by spherical caps of radius $\varphi/2$.

The bound in Theorem 1.1 yields
$$
A(3,\varphi)\le \frac{2\pi}{\Delta(\varphi)}+2,
$$
where
$$
\Delta(\varphi)=3\arccos{\left(\frac{\cos{\varphi}}{1+\cos{\varphi}}\right)}-\pi
$$
is the area of a spherical regular triangle with side length $\varphi$.

The bound is tight for $N=3,4,6$ and 12. So for these $N$ it gives a solution of the Tammes problem. It is also tight asymptotically.
However, is not tight for any other cases. 

\subsection{Coxeter's bound}
In 1963 Coxeter  \cite{Cox} proposed an extension of Fejes T\'oth's bound for all dimensions. His bound was based on the conjecture that in $n$-dimensional spherical space equal size balls cannot be packed denser than the density of $n+1$ mutually touching balls of the same size with respect to the simplex spanned by the centers of the balls. This conjecture has been stated by Fejes T\'oth for the 3-dimensional case in \cite{FeT53} and  for all dimensions in  \cite{FeT59, FeT60}. Assuming the correctness of the conjecture Coxeter calculated the upper bounds 26, 48, 85, 146, and 244 for the kissing numbers $k(n)=A(n,\pi/3)$ for $ n=4, 5, 6, 7,$ and 8, respectively. The conjecture, which was finally confirmed by B\"or\"oczky \cite{Bor1} in 1978 also yields that
$$
A(4,\pi/5)=120.
$$

\begin{theorem}[B\"or\"oczky \cite{Bor1} and Coxeter  \cite{Cox}]
$$
A(n,\varphi)\le 2F_{n-1}(\alpha)/F_n(\alpha),	
$$
where
$$
\sec{2\alpha}=\sec{\varphi}+n-2,
$$
and the function $F$ is defined recursively by
$$
F_{n+1}(\alpha)=\frac{2}{\pi}\int\limits_{{arcsec}(n)/2}^\alpha{F_{n-1}(\beta)\,d\theta}, \; \sec{2\beta}=\sec{2\theta}-2,
$$
with the initial conditions $F_0(\alpha)=F_1(\alpha)=1$.
\end{theorem}

\section{The problem on maximizing the minimum distance between antipodal pairs of points on the sphere}

L. Fejes T\'oth \cite{FeT65} considered Tammes' problem for antipodal sets on ${\mathbb S}^2$ {\it{i.e.}} for sets $X$ that are invariant under the antipodal mapping $A: {\mathbb S}^d\to {\mathbb S}^d$, where \; $A(x)=-x$.
Let
$$
a_M:=\max\limits_{X=-X\subset{\Bbb S}^2}{\{\psi(X)\}}, \, \mbox{ for } \;  |X|=2M.
$$
For a given $M$, Fejes T\'oth's problem for antipodal sets is to find all configurations
$$
X=\{x_1,-x_1,\ldots,x_M,-x_M\}
$$
on ${\mathbb S}^2$ such that $\psi(X)=a_M$.

This problem is presently solved only for $M\le 7$. It is clear that $a_M\le d_{2M}$. Therefore, if $\psi(X)=d_{2M}$, $|X|=2M$, and $X$ is antipodal then $a_M=d_{2M}$. Thus, for $M=3$ and $M=6$ we have this equality. The following theorem is the main result of \cite{FeT65}.

\begin{theorem}[L. Fejes T\'oth, \cite{FeT65}]
Let $P_M\subset{\Bbb S}^2$ be a maximal set for the Fejes T\'oth  problem for antipodal configurations, i. e.  $\psi(P_M)=a_M$. Then
		\begin{enumerate}
		\item $P_2$  is the set of vertices of a square on the equator, $a_2=90^\circ$;
		\item   $P_3$ is the set of vertices of a regular octahedron, $a_3=90^\circ$;
		\item   $P_4$ is the set of vertices of a cube, $a_4=\arccos{(1/3)}$;
		\item    $P_5$ consists of five pairs of antipodal vertices of a regular icosahedron,  $a_5=\arccos{(1/\sqrt{5})}$.
		\item    $P_6$ is the set of vertices of a regular icosahedron, $a_6=\arccos{(1/\sqrt{5}})$.
	\end{enumerate}
\end{theorem}

In our paper with Tarasov \cite{MTEP} we gave an alternative proof of this theorem. In \cite{MT2013} we found the list of all irreducible contact graphs with $N$ vertices on the sphere ${\mathbb S}^2$, where $6\le N \le 11$. Since the contact graph of $P_M$ is irreducible the theorem (for $M<6$) follows from this list.

Fejes T\'oth conjectured that the solution of the problem for seven pairs of antipodal points consists of the vertices of a rhombic dodecahedron (see the second edition of \cite{FeT}, page 210). This was proved by Cohn and  Woo \cite{CW12} as a consequence of a more general theorem.
\medskip

\section{Problems on the maximum contact number of packings on the sphere}

In \cite{FeT86} (pages 86 and 87) Fejes T\'oth raised three problems abot the number of touching pairs in a packing of congruent circles on the sphere.

Consider a packing $P$ of  ${\mathbb S}^2$ by $N$ circles $c_1,\ldots,c_N$ of diameter $d$.
In the packing  $P$ let $c_i$ be touched by $k_i$ circles. The first problem is to {\it find the maximum  number of points of contact:}
$$
\ko_N(d):=\max\limits_{P: |P|=N}\frac{k_1+\ldots+k_N}{2}.
$$
In other words, $\ko_N(d)$ is the maximum number of touching pairs in a packing of $N$ spherical caps of diameter $d$.

It is clear that if  $d=d_N$, then $\ko_N(d)$ is realized by the solution of the Tammes problem. It seems that the case $d<d_N$ is not well considered. There is only one paper in this direction \cite{FTTT}, where this problem is considered for $N = 12$ and $d = 60^\circ$. There, it is proved that
$$\ko_{12}(60^\circ)=24.$$

Fejes T\'oth also proposed the problem of {\it finding the maximum
$$\bar K(d):=\max\limits_N{\frac{\ko_N(d)}{N}}$$
of the average number of points of contacts over all packings of circles of diameter $d$.}

The third problem is: {\it For a given $N$, find the maximum kissing number $\ko_N$ over all packings of equal circles,} {\it{i.e.}}, find
$$
\ko_N:=\max\limits_{d\le d_N}{\ko_N(d)}.
$$

Let $X$ be the set of centers of a packing of congruent circles on ${\Bbb S}^2$. Denote by $e(X)$ the number of edges of the contact graph $\cg(X)$. It is easy to see that
$$
\ko_N=\max\limits_{X\in {\Bbb S}^2, |X|=N}{e(X)}.
$$
This number is currently known only for $N\le12$ and $N=24,48,60,120.$

Denote by $\kappa(d)$ the kissing number of the spherical cap with diameter $d$ in ${\Bbb S}^2$, {\it{i.e.}} it is the maximum number of non-overlapping circles of diameter $d$ that can touch a circle of the same diameter. Note that if $d\le\arccos(1/\sqrt{5})$, then $\kappa(d)=5$.

We say that a packing of $N$ spherical caps with diameter $d$ is {\it maximal} if
$$\ko_N(d)=N\kappa(d)/2.$$
The following theorem has been proved by Robinson \cite{Rob69} and  Fejes T\'oth \cite{FeT69}.

\begin{theorem}[Robinson \cite{Rob69}, Fejes T\'oth \cite{FeT69}] A maximal packing of $N$ equal spherical caps exists only if $N=2,\,3,\,4,\,6,\,8,\,9,\,12,\,24,\,48,\,60$ or $120$.
\end{theorem}

This theorem  implies
\begin{cor}
 $\ko_2=1,\, \ko_3=3,\, \ko_4=6,\, \ko_6=12,\, \ko_8=16, \, \ko_9=18$ and for $N=12,\,24,\,48,\,60$  or $120$ we have $\ko_N=5N/2.$
 \end{cor}

In our paper \cite{MTEP} we considered $\ko_N$ for $N<12$. In particular, we proved that

 \begin{theorem} [Musin and Tarasov \cite{MTEP}] \label{t32}
 $\ko_5=8$, $\ko_7=12$,  $\ko_{10}=21,$ and  $\ko_{11}=25$.
\end{theorem}

Note that $\ko_5$ is attained by the set of vertices of a square pyramid. For $N=7$ and $N=11$,  $\ko_N(d)$ achieves its maximum on optimal configurations for Tammes' problem. However, the arrangement realizing the optimal value $\ko_{10}$ is obtained by removing from the set of vertices of a regular icosahedron two adjacent vertices. In this case the contact graph $\cg(X)$ is not irreducible.

Our proof of Theorem \ref{t32} in \cite{MTEP} is based on two lemmas.

\begin{lemma}
Let $X$ be a finite set on the sphere ${\mathbb S}^2$. If every face of the contact graph $\cg(X)$
is either a triangle or a quadrilateral, then this graph is irreducible.
\end{lemma}

 \begin{lemma}  Let $X$ be a finite set on the sphere ${\mathbb S}^2$, where
 $|X| = N$ and $N > 6.$ Suppose that $e(X) \ge 3N - 8.$ Then the
contact graph $\cg(X)$ is irreducible.
 \end{lemma}
Using these lemmas, Theorem \ref{t32} follows by checking the list of irreducible contact graphs for $N\le11$ \cite{MT2013}.

\section{The Fejes T\'oth -- Sachs problem on the one--sided kissing numbers}

Let $H$ be a closed half-space of ${\Bbb R}^{n}$. Suppose $S$ is a unit sphere in $H$ that touches the bounding hyperplane of $H$. The {\it one--sided kissing number} $B(n)$ is the maximal number of unit non--overlapping spheres in $H$ that can touch $S$.


The problem of finding $B(3)$ was raised by Fejes T\'oth and Sachs in 1976 \cite{FeTS76} in another context.
K. Bezdek and Brass \cite{BB} studied the problem in a more
general setting and they introduced the term ``one-sided Hadwiger
number", which in the case of a ball is the same as the one-sided kissing
number. The term ``one-sided kissing number" has been introduced by K. Bezdek \cite{KB}.

Clearly, $B(2)=4$.   The Fejes T\'oth -- Sachs problem in three dimensions was solved by G. Fejes T\'oth \cite{GFT}. He proved that $B(3)=9$ (see also Sachs \cite{Sachs} and A. Bezdek and K. Bezdek \cite{AKB} for other proofs). Finally, Kert\'esz \cite{Ker} proved that the maximal one--sided kissing arrangement is unique up to isometry.

The first upper bound for $B(4)$ was given by Szab\'o \cite{Szabo}. He used the Odlyzko--Sloane bound $k(4)\le 25$ for the kissing number of the four-dimensional ball to show that $B(4) \le 20$. Next K. Bezdek \cite{KB}, based on the result that $k(4)=24$ \cite{Mus,Mus2}, lowered the bound to $B(4)\le 19$.

In  \cite{Mus3} I proved that $B(4)=18$. This proof relies on the extension of Delsarte's method that was developed in \cite{Mus2}. However, technically the proof is more complicated than the proof of the fact that $k(4)=24$.  An alternate proof was given in \cite{BV2} using semidefinite programming. The problem of uniqueness of the maximal one--sided kissing arrangement in four dimensions is still open.

 In \cite{Mus3} I conjectured that $B(5)=32, \; B(8)=183$ and  $B(24)=144855$. This conjecture for $n=8$ was proved by  Bachoc anf Vallentin \cite{BV2}. In \cite{BM} and  \cite{Mus4} we proposed several upper bounds on $B(n)$. However, all these bounds were improved in \cite{BV2}.


It is clear that there are some relations between kissing numbers and one--sided kissing numbers.
Look at these nice equalities:

$$ n=2, \quad 4=B(2)=\frac{k(1)+k(2)}{2}=\frac{2+6}{2};$$
$$ n=3, \quad 9=B(3)=\frac{k(2)+k(3)}{2}=\frac{6+12}{2};$$
$$ n=4, \quad 18=B(4)=\frac{k(3)+k(4)}{2}=\frac{12+24}{2}.$$

We do not know whether the equality
$$B(n)=\bar K(n):=\frac{k(n-1)+k(n)}{2}$$
holds for all $n$.  However, there are  reasons to believe that $B(n)=\bar K(n)$ for
$ n=5, 8$ and $24$.  We propose a weaker conjecture, namely, that the equality $B(n)=\bar K(n)$ holds asymptotically:
\medskip

\noindent{\bf Conjecture.} We have
$$ \lim_{n\to\infty}{\frac{B(n)}{\bar K(n)}}=1.$$

\section{The work of Fejes T\'oth on the isoperimetric problem for polyhedra and their extensions}

\subsection{Isoperimetric problem for polyhedra}

The isoperimetric problem  in space can be formulated as follows:  {\it Find a convex body of given surface area $F$ which contains the largest volume $V$.}

The famous isoperimetric inequality states
$$
F^3\ge 36\pi V^2.
$$

For any solid $P$ consider the {\it Isoperimetric
Quotient}
$$
\iq(P)= 36\pi\frac{V^2}{F^3},
$$
 a term introduced by P\'olya in \cite[Chap. 10, Problem 43]{Pol}.  The isoperimetric inequality implies that $\iq(P)\le 1$ and the equality  holds only if $P$ is a sphere.

The isoperimetric problem for polyhedra was first considered by Lhuilier (1782), see \cite{Lhu}, and Steiner (1842), see \cite{Steiner}. 
Steiner stated the following conjecture. 

\medskip

\noindent{\bf Steiner's conjecture \cite{Steiner}.}  {\it Each of the five Platonic solids is the best (i.e. with the highest $\iq$) among all isomorphic polyhedra.}

\medskip

This problem is still open for the icosahedron.

\medskip

Consider the isoperimetric problem for polyhedra with given number of faces $f$.
Actually, this problem is currently solved only for $f\le 7$ and $f=12$.  However, the first theorem on this problem was discovered in the 19th century.

\begin{theorem} [Lindel\"of \cite{Lin} and Minkowski \cite{Min}] {Of all convex polyhedra with the same number of faces, a polyhedron with the highest $\iq$ is circumscribed about a sphere which touches each face in its centroid.}
\end{theorem}

Note that IQ(tetrahedron)$\approx 0.302$, IQ(cube)$\approx 0.524$, IQ(octahedron)$\approx 0.605$, IQ(dodecahedron)$\approx 0.755$, and IQ(icosahedron)$\approx 0.829$  \cite[Chap. 10, p. 189]{Pol}. 

In fact, there are simple polyhedra $F_{12}$ and $C_{36}$ with eight and 20 faces that have greater IQ than, respectively, the regular octahedron and icosahedron. Goldberg \cite{Gold} computed that $\iq$(octahedron)$<\iq(F_{12})\approx 0.628$ and $\iq$(icosahedron)$<\iq(C_{36})\approx 0.848$. 

Goldberg proved that the regular dodecahedron is the best polyhedron
with 12 facets and stated the following conjecture.

\medskip

\noindent{\bf Goldberg's conjecture \cite{Gold}.} {\it If a polyhedron $P$ with $f\ne 11,13$ faces and $v$ vertices has the greatest $\iq$, then $P$ is simple and its faces are $\lfloor 6-24/(v+4)\rfloor$--gons or  $\lfloor 7-24/(v+4)\rfloor$--gons.}

\medskip

Note that according to Goldberg's conjecture if $v\ge 20$, then the faces of a best polyhedron can be only pentagons and hexagons, in other words $P$ is a {\it fullerene} (see \cite{D2G}).

Let $P$ be a convex polyhedron with $f$ faces. In \cite{Gold}  Goldberg proposed the following ineaqulity:
$$
\frac{F^3}{V^2}\ge 54\,(f-2)\tan{\omega_f}\,(4\sin^2{\omega_f}-1),
\;  \; \omega_f:=\frac{\pi f}{6f-12},
$$
or equivalently
$$
\iq(P)\le \frac{2\pi\cot{\omega_f}} {3(f-2)(4\sin^2{\omega_f}-1)}, \eqno (5.1)
$$
where the equality holds only if $f=4$ (regular tetrahedron), $f=6$ (cube) or $f=12$ (regular dodecahedron).

This inequality was independently found by Fejes T\'oth \cite{FeT44} (see Fejes T\'oth's books \cite{FeT, FeTbook2} and Florian's survey \cite{Florian} for references and historical remarks).  However, both proofs  contained a gap, namely the proof of the convexity of a certain function of two variables. A first rigorous proof of $(5.1)$ was given by Fejes T\'oth in the paper \cite{FeT48}. Finally, Florian \cite{Flor56} filled the gap in the previous proof by establishing the convexity of the respective function.

Two conjectures of  Fejes T\'oth on isoperimetric inequalities are still open.  Let $P$ be a convex polyhedron with $v$ vertices. The first conjecture states that
$$
\frac{F^3}{V^2}\ge \frac{27\sqrt{3}}{2}(v-2)(3\tan^2{\omega_v}-1).
$$

Let $P$ be a convex polyhedron with $f$ faces, $v$ vertices and $e$ edges. The second conjecture of Fejes T\'oth states that
$$
\frac{F^3}{V^2}\ge 9\,e\sin{\frac{2\pi}{p}}\left(\tan^2{\frac{\pi}{p}}\, \tan^2{\frac{\pi}{q}}-1\right),
$$
where $p:=2e/f$ and $q:=2e/v$.

Note that the validity of any of these conjectures would yield a proof of the still open conjecture of Steiner concerning the isoperimetric property of the icosahedron.

\subsection{The Goldberg --  Fejes T\'oth inequality}

Here we consider a ``dual'' version of the Goldberg -- Fejes T\'oth inequality $(5.1)$  

Let $P$ be a convex polyhedron with $f$ faces. Then Euler's formula implies
$$
v\le 2f-4,  \eqno (5.2)
$$
where $v$ is the number of vertices of $P$. The equality holds only if $P$ is a simple polyhedron.

Suppose $P$ is a polyhedron with highest $\iq$ and fixed $f$. Then by the Lindel\"of -- Minkowski theorem there is a sphere $S$ that touches each face $\Gamma_i$ of $P$ in its centroid $x_i$.
Thus, the set   $X:=\{x_1,\ldots,x_f\} $ is a subset of $S$.

Without loss of generality it can be assumed that $S$ is of radius $r=1$.  Since all faces $\Gamma_i$ touch the unit sphere $S$,  we have 
$$
V=\frac{1}{3}F. 
$$
Therefore, for a given $f$, $P$ has the highest $\iq$ if and only if $F=\ar(P)$ achieves its minimum.

Let $O$ be the center of the sphere $S$.   Consider the central projection $g:P\to S$, i.e. for any point $A\in P$ the projection $g(A)$ is the intersection of the line $OA$ with $S$. It is clear that  $g(x_i)=x_i$.

Let  $p_1,\ldots,p_v$ be vertices of $P$. Denote by $Q$ the projection of this vertex set, i.e. $Q:=\{q_1,\ldots,q_v\}$, where $q_i:=g(p_i)\in S$.

The set $Q$ coincides with the set of vertices $\{v_i\}$ of the Voronoi diagram $\vd(X)$ of  $X$  in $S$. Equivalently, if $G_{i1},\ldots,G_{im}$ are faces of $\vd(X)$ with a common vertex $v_i$, then $q_i$ is the circumcenter of the Delaunay cell $D_i$ with vertices $x_{i1},\ldots,x_{im}$ in $S$. It immediately follows from the fact that $|v_ix_{ij}|$ as well as $|q_ix_{ij}|$ does not depend on $j$.  Indeed, we have
$$
|v_ix_{ij}|^2=|Ov_i|^2-|Ox_{ij}|^2, \; \mbox{ where } \;  |Ox_{ij}|=1  \; \mbox{ for all } \; j=1,\ldots, m.
$$

Let $G_i:=g(\Gamma_i)$, $i=1,\ldots,f$. Then $G_i$ are Voronoi cells of $\vd(X)$. Let $D_1,\ldots,D_v$ be the cells of the Delaunay tessellation $\dt(X)$ in $S$. Then 
$$
\ar(G_1)+\ldots+\ar(G_f)=\ar(D_1)+\ldots+\ar(D_v)=4\pi. \eqno (5.3)
$$

We have 
$$
F=\ar(\Gamma_1)+\ldots+\ar(\Gamma_f)=\ar(g^{-1}(G_1))+\ldots+\ar(g^{-1}(G_f)). 
$$
It is easy to see, that the equality in  the Goldberg --  Fejes T\'oth inequality (see \cite{Gold} and  \cite[Sect. V.4]{FeT}) holds only if $v=2f-4$ and all $G_i$ are congruent regular polygons. Equivalently, that means the equality holds only if all $D_k$ are congruent regular triangles.

Denote $\tilde D_i=g^{-1}(D_i)$. Then 
$$
F=\ar(\tilde D_1)+\ldots+\ar(\tilde D_v). \eqno (5.4)
$$
If $v=2f-4$ and all $D_i$ are  congruent triangles, then (5.3) yields
$$
\ar(D_i)=\frac{2\pi}{f-2}, \quad i=1,\ldots,v. 
$$

Let $T$ be a regular spherical triangle in $S$ with $\ar(T)=t$.
Denote $$\rho(t):=\ar(\tilde T), \quad \tilde T:=g^{-1}(T)\subset P_T.$$ 
Thus, we have the following inequality that is equivalent to $(5.1)$.
$$
\iq(P)\le \frac{\tau}{\rho(\tau)}, \; \; \tau=\frac{2\pi}{f-2}. \eqno (5.5)
$$

\subsection{The Goldberg --  Fejes T\'oth inequality for higher dimensions}

Here we consider a $d$--dimensional analog of the inequality $(5.5)$.

Let $P$ be a convex polyhedron  in ${\Bbb R}^d$ with $v$ vertices and $n$ facets, i.e. $v=f_0(P)$ and  $n=f_{d-1}(P)$. Then
McMullen's Upper Bound Theorem \cite[p. 254]{Ziegler}  yields the extension
$$
v\le h_d(n):= {n- \lceil d/2\rceil \choose \lfloor d/2 \rfloor}+{n- \lfloor d/2\rfloor -1 \choose \lceil d/2 \rceil -1}. \eqno (5.6)
$$
of the inequality $(5.2)$ for all dimensions.

Define
$$
\iq(P):=d^{d-1}\,\Omega_d\,\frac{V^{d-1}}{F^d}, \; \; \Omega_d:=\ar({\Bbb S}^{d-1}).
$$

The Lindel\"of -- Minkowski theorem holds for all dimensions (see \cite[p. 274]{Hadw}).  So we have that  $\iq(P)$ achieves its maximum only if $P$ is circumscribed about a sphere $S$.  As above we assume that $S$ is a unit sphere.

It is easy to see that all definitions from Subsection 5.2 can be extended to all dimensions. Therefore, equality  $(5.4)$ holds also  for a $d$-dimensional polyhedron $P$. An analog of $(5.3)$ is the following equality:
$$
\ar(D_1)+\ldots+\ar(D_v)=\Omega_d. \eqno (5.7)
$$

Let $D$ be a regular spherical simplex in $S$ with $\ar(D)=t$.
Denote $$\rho_d(t):=\ar(\tilde D).$$ Our conjecture is that the following extension of the Goldberg --  Fejes T\'oth inequality $(5.5)$ holds for all dimensions:
$$
\iq(P)\le \frac{{\Omega_d}}{v\,\rho_d({\Omega_d}/{v})}.  \eqno (5.8)
$$
Since $v\le h_d(n)$, in particular we have 
$$
\iq(P)\le \frac{\tau}{\rho_d(\tau)}, \; \; \tau=\frac{\Omega_d}{h_d(n)}. \eqno (5.9)
$$

Perhaps, $(5.8)$ can be proved by the same way as L\'aszl\'o  Fejes T\'oth proved $(5.1)\equiv (5.5)$ in \cite{FeT48} and  \cite[Sect. V.4]{FeT}. Actually, for all dimensions there are extensions of $(5.2)-(5.4)$.  We think that the most complicated step here is to prove that $\iq(P)$ cannot exceed $\iq$ of a such polyhedron with $n$ facets  that all its $D_i$ are congruent  regular spherical simplices.

\medskip

\medskip

\medskip

\noindent{\bf Acknowledgment.} I wish to thank  G\'abor Fejes T\'oth for useful comments, references,  and thoughtful editing of my paper.

\medskip

\medskip

\medskip

\medskip

 O. R. Musin,  University of Texas Rio Grande Valley, School of Mathematical and Statistical Sciences, One West University Boulevard, Brownsville, TX, 78520.

 {\it E-mail address:} oleg.musin@utrgv.edu

\end{document}